\def\titlep{Algebra with Indefinite Involution and Its Representation
in Krein Space}
\documentclass[11pt]{article}
\usepackage{graphicx,ifthen}
\usepackage{amssymb}
\usepackage{amsmath}

\font\germ=eufm10 at12pt

\def\goth#1{\hbox{\germ#1}}

\newcommand{\sdag}{\scriptsize \dag}
\newcommand{\cpd}{{\hbox{\math o}}}
\font\math=msbm10 scaled 1200

\newcommand{\qed}{\hbox{\rule[-2pt]{3pt}{6pt}}}
\newcommand{\qedh}{\hfill\qed \\}

\newtheorem{Thm}{Theorem}

\newtheorem{ex}{Example}
\newtheorem{defi}{Definition}
\newtheorem{lem}{Lemma}

\newtheorem{prop}{Proposition}

\newtheorem{cor}{Corollary}

\def\cal#1{\mathcal #1}
\def\con{{\cal O}_{N}}

\def\edot{=1,\ldots,N}
\def\pr{{\it Proof.}\quad}

\def\co#1{{\cal O}_{#1}}
\def\ltn{l_{2}({\bf N})}

\def\ggl{{\goth g}{\goth l}}

\def\gu{{\goth u}}
\def\bfsnl{{\rm BFS}_{N}(\Lambda)}

\addtocounter{footnote}{1}
\def\sftt#1{
\setcounter{equation}{0}
\addtocounter{footnote}{1}
\section{#1}
}

\def\ssfr#1{\subsection*{#1}}

\begin{document}
%
%
\def\autherp{Katsunori Kawamura}
\def\emailp{{\it E-mail: kawamura@kurims.kyoto-u.ac.jp}}
\def\addressp{{\it College of Science and Engineering Ritsumeikan University,}\\
{\it 1-1-1 Noji Higashi, Kusatsu, Shiga 525-8577, Japan.}
}
\def\bfsnl{{\rm BFS}_{N}(\Lambda)}
\def\scm#1{S({\bf C}^{N})^{\otimes #1}}
\def\mqb{\{(M_{i},q_{i},B_{i})\}_{i=1}^{N}}
\newcommand{\mline}{\noindent
\thicklines
\setlength{\unitlength}{.1mm}
\begin{picture}(1000,5)
\put(0,0){\line(1,0){1250}}
\end{picture}
\par
 }
\def\sd#1{#1^{\sdag}}
\def\pca{pseudo-Cuntz algebra}
\def\emailp{e-mail: kawamura@kurims.kyoto-u.ac.jp.}
\def\addressp{College of Science and Engineering Ritsumeikan University,\\
1-1-1 Noji Higashi, Kusatsu, Shiga 525-8577, Japan}
%
%
\setcounter{section}{0}
\setcounter{footnote}{0}
\setcounter{page}{1}
\pagestyle{plain}
\title{\titlep}
\author{\autherp\thanks{\emailp}
\\
\addressp}

%
%


\maketitle
%
%
\begin{abstract}
It is often inevitable to introduce an indefinite-metric space
in quantum field theory,  for example,
which is explained for the sake of 
the manifestly covariant quantization of the electromagnetic field. 
We show two more evident mathematical reasons why such indefinite metric 
appears.
The first idea is the replacement of involution on an algebra.
For an algebra ${\cal A}$ with an involution $\sdag$
such that a representation of the involutive algebra $({\cal A},\sdag)$  
brings an indefinite-metric space,
we replace the involution $\sdag$ with a new one $*$ on ${\cal A}$
such that $({\cal A},*)$ is a well-known involutive algebra
acting on a representation space with positive definite metric.
This explains that non-isomorphic two involutive algebras 
are transformed each other by
the replacement of involution.
The second is that a covariant (Hilbert space) representation 
$({\cal H},\pi,U)$ of an
involutive dynamical system $(({\cal A},*),{\bf Z}_{2},\alpha)$
brings a Krein space representation 
of the algebra ${\cal A}$ with the replaced involution.
For example,
we show representations of abnormal CCRs, CARs and pseudo-Cuntz algebras 
arising from those of standard CCRs, CARs and Cuntz algebras.
\end{abstract}

\sftt{Introduction}
\label{section:first}
We have studied representations of operator algebras 
\cite{AK01,AK05,IWF01,IMQ01,RBS01,TS01}.
Most of the representation spaces are complex vector 
spaces with positive definite metric. 
However, it is often inevitable to introduce an indefinite-metric space
in quantum field theory \cite{Nagy,Nakanishi}.
In this paper, 
we clarify the mechanism why such indefinite metric appears
and show a systematic treatment of such indefinite-metric representation
from a standpoint of involutions on algebras.

%
%
\subsection{Involution and vector space with indefinite metric}
\label{subsection:firstone}
We explain terminology here.
In this paper, any algebra means an algebra over ${\bf C}$.
A map $\varphi$ on ${\cal A}$ is called an {\it involution} 
on ${\cal A}$ if $\varphi$ is a conjugate linear map which satisfies 
$\varphi(xy)=\varphi(y)\varphi(x)$ for each $x,y\in {\cal A}$ 
and $\varphi^{2}=id$.
For convenience,
we write $x^{\varphi}$ instead of $\varphi(x)$ for $x\in {\cal A}$.
For the notation of involution, $*,\sdag$ and $\#$ are often used
\cite{NakagamiTomita,Pauli,Tomita}.
In physics, for an operator $T$, the operator $T^{\sdag}$ is often called
the {\it hermite conjugate} of $T$, which is the image of 
an involution $\sdag$ of $T$.
Remark that a different notion of ``involution" is used in 
supersymmetry ($\S$ 5.1.1 in \cite{Thaller}).

The terminology of ``$*$-algebra" is not suitable
to treat two different involutions on an algebra at once.
Hence we use a terminology, involutive algebra instead of it
according to Chapter 1 of \cite{Bourbaki}.
A pairing $({\cal A},*)$ is an {\it involutive algebra}
if $*$ is an involution on an algebra ${\cal A}$.
Of course, a $*$-algebra ${\cal A}$ is an involutive algebra $({\cal A},*)$.
If $x\in {\cal A}$ satisfies $x^{*}=x$,
then $x$ is called {\it $*$-self-adjoint}.
An involutive algebra $({\cal A},*)$ is called a {\it Banach involutive algebra} 
if ${\cal A}$ is a Banach algebra and $\|x^{*}\|=\|x\|$ for each 
$x\in {\cal A}$.
A Banach involutive algebra $({\cal A},*)$ is a $C^{*}$-algebra
if $\|x^{*}x\|=\|x\|^{2}$ for each $x\in {\cal A}$.
For involutive algebras $({\cal A},*)$ and $({\cal B},\sdag)$, 
a homomorphism $f$ from ${\cal A}$ to ${\cal B}$ is {\it involutive}
if $f\circ *=\sdag\circ f$.
An automorphism $\alpha$ of $({\cal A},*)$
is {\it involutive} if $*\circ \alpha=\alpha\circ *$.
An involution $\sdag$ on ${\cal A}$ is {\it equivalent} to $*$ if
there exists an involutive isomorphism from
$({\cal A},\sdag)$ to $({\cal A},*)$.
A subalgebra ${\cal A}_{0}$ of $({\cal A},*)$ is 
{\it involutive} (or {\it self-adjoint})
if $\{x^{*}:x\in {\cal A}_{0}\}\subset {\cal A}_{0}$.

An involution is one of most important structures on operator algebras
\cite{Naimark}. Especially, the 
C$^{*}$-condition is a nice characterization of a special involution 
with respect to the norm.
On the other hand,
properties of the involution on the algebra of 
field operators in quantum field theory is not well-known.
For example, 
the involution on the algebra ${\cal A}$ of field operators 
in quantum electrodynamics satisfies
neither the C$^{*}$-condition nor the positivity of 
the spectrum of the operator $I+x^{*}x$ for $x\in {\cal A}$.

On a complex vector space $V$, a map $(\cdot|\cdot)$ 
from $V\times V$ to ${\bf C}$
is called a {\it hermitian form} on $V$ 
if $(\cdot|\cdot)$ is sesquilinear and $\overline{(v|w)}=(w|v)$ 
for each $v,w\in V$.
We assume that $(x|cy)=c(x|y)$ for $x,y\in V$ and $c\in {\bf C}$
in this paper.
In the theory of indefinite-metric space and quantum field theory,
such a hermitian form is called by an {\it inner product} or {\it metric}.
A hermitian form $(\cdot|\cdot)$ is {\it indefinite}
if there exist $v,w\in V$ such that $(w|w)<0<(v|v)$.
Such pair $(V,(\cdot|\cdot))$ is called an {\it indefinite-metric space}
or {\it indefinite-inner product space}.

For any operator $T$ on a hermitian vector space $(V,(\cdot|\cdot))$,
if $(\cdot|\cdot)$ is nondegenerate, 
then there exists unique operator $T^{\star}$ on $V$
such that $(T^{\star}v|w)=(v|Tw)$ for each $v,w\in V$.
The involution $\star$ is called the {\it involution
associated with the hermitian form $(\cdot|\cdot)$}. 
With respect to this $\star$,
the algebra ${\rm End}V$
of all linear operators on $V$ is an involutive algebra $({\rm End}V,\star)$.
A pairing $(V,\pi)$ is
a {\it (involutive) representation} of an involutive algebra $({\cal A},\sdag)$ 
if $V$ is a vector space with a nondegenerate hermitian form $(\cdot|\cdot)$
and $\pi$ is an involutive  
homomorphism from ${\cal A}$ to $({\rm End}V,\star)$.

%
%
\subsection{Abnormal commutation relations}
\label{subsection:firsttwo}
According to the preface in \cite{Bognar},
the theory of indefinite-metric space 
has two origins which are relatively independent. 
One is quantum field theory \cite{Dirac,Pauli}
and other is functional analysis \cite{Pontrjagin,Sobolev}.
A difference of their styles is the order of logic of the theory.
In the former, 
the algebra ${\cal A}$ of field operators appears in first and
an indefinite-metric space appears as a representation space of ${\cal A}$. 
On the other hand, the later, 
an indefinite-metric space is given in first and 
the theory of linear operators on it is discussed \cite{MS1980,MPS1990,MPS1992,Strocchi}.

We show that an algebras with a certain involution 
brings an indefinite-metric space as its involutive representation.
In order to explain such algebras, we demonstrate by three examples.

For three families $\{a,a^{\sdag}\}, \{f,f^{\sdag}\}$ and 
$\{s_{1},s_{2},s_{1}^{\sdag},s_{2}^{\sdag}\}$ with an involution $\sdag$,
consider the following relations:
%
%
\begin{equation}
\label{eqn:abnormalboson}
aa^{\sdag}-a^{\sdag}a=-I, \quad aa-aa=a^{\sdag}a^{\sdag}-a^{\sdag}a^{\sdag}=0,
\end{equation}
%
%
\begin{equation}
\label{eqn:abnormalfermion}
ff^{\sdag}+f^{\sdag}f=-I,\quad ff+ff
=f^{\sdag}f^{\sdag}+f^{\sdag}f^{\sdag}=0,
\end{equation}
%
%
\begin{equation}
\label{eqn:pseudo}
s_{i}^{\sdag}s_{j}=(-1)^{i-1}\delta_{ij}I\quad(i,j=1,2),
\quad s_{1}s_{1}^{\sdag}-s_{2}s_{2}^{\sdag}=I
\end{equation}
where $I$ denotes the unit in each case.
Relations (\ref{eqn:abnormalboson}),
(\ref{eqn:abnormalfermion}) and (\ref{eqn:pseudo})
are called 
the {\it abnormal canonical commutation relations},
the {\it abnormal canonical anti-commutation relations},
and the {\it pseudo-Cuntz relations}, respectively
\cite{AK05,Nakanishi}.
Define involutive algebras 
${\cal A}_{\bar{B}}$, ${\cal A}_{\bar{F}}$ and $\co{1,1}^{(0)}$ 
generated by them, respectively.
We consider their involutive representations as follows.

In Section  3 of \cite{Nakanishi},
it was explained that 
involutive representations of ${\cal A}_{\bar{B}}$ and ${\cal A}_{\bar{F}}$ 
cause indefinite-metric spaces as the involutive representation spaces of them.
Assume that $(\cdot|\cdot)$ is a nondegenerate hermitian form on $V$
and ${\cal A}_{\bar{B}}$ is involutively 
represented on $(V,(\cdot|\cdot))$.
If there exists a vector $\Omega\in V$ such that  $(\Omega|\Omega)>0$ and $a\Omega=0$,
then $(a^{\sdag}\Omega|a^{\sdag}\Omega)=-(\Omega|\Omega)<0$.
Hence (\ref{eqn:abnormalboson}) brings
an indefinite-metric representation in this case.
Because the algebra generated by (\ref{eqn:abnormalboson})
is involutively isomorphic to that of canonical commutation relations,
the reason why an indefinite-metric space appears
may be considered as the choice of the vacuum vector $\Omega$.
On the other hand,
we can verify that any unital involutive representation of ${\cal A}_{\bar{F}}$
must be an indefinite-metric space.
In \cite{IMQ01}, we introduced
$\eta$-CCRs and $\eta$-CARs which are generalization 
of (\ref{eqn:abnormalboson}) and
(\ref{eqn:abnormalfermion}), and  their representations on Krein spaces
by modifying Fock representations of standard CCRs and CARs.

In \cite{AK05}, we introduced  $\co{1,1}^{(0)}$
which is called the {\it pseudo-Cuntz algebra},
in order to construct representations of 
the Faddeev-Popov (=FP) (anti) ghost fields in string theory.
It is understood that there is no C$^{*}$-algebra
which contains $\co{1,1}^{(0)}$ as an involutive subalgebra.
We construct an involutive representation with indefinite metric
of $\co{1,1}^{(0)}$ as follows.
Consider a representation
of $(V,\pi)$ of $\co{1,1}^{(0)}$ with a cyclic
vector $\Omega$ 
such that
\[\pi(s_{1})\Omega=\Omega.\]
By (\ref{eqn:pseudo}),
we see that $\pi(s_{2}^{\sdag})\Omega=0$.
Hence the cyclic representation space $V$ of $\co{1,1}^{(0)}$ 
is the linear span of the family $\{\Omega, \pi(s_{j_{1}}\cdots s_{j_{k}}s_{2})
\Omega:j_{1},\ldots,j_{k}=1,2\mbox{ for }k\geq 1\}$
of vectors.
Define the hermitian form on $V$ by
\[(e_{J}|e_{K})=(-1)^{n_{2}(J)}\delta_{JK}\]
where $e_{J}\equiv 
\pi(s_{j_{1}}\cdots s_{j_{k}})\Omega$ and 
$n_{2}(J)\equiv \sum_{i=1}^{k}(j_{i}-1)$
for $J=(j_{1},\ldots,j_{k})$.
Then $(\cdot|\cdot)$ is nondegenerate on $V$
and $(\Omega|\Omega)=1$.
We see that $\co{1,1}^{(0)}$ acts on $(V,(\cdot|\cdot))$ involutively.
Hence $((V,(\cdot|\cdot)),\pi)$ is 
an involutive representation with indefinite metric
of $\co{1,1}^{(0)}$.
%
%
\subsection{Replacement of involution}
\label{subsection:firstthree}
We replace the involution $\sdag$ 
on ${\cal A}_{\bar{B}},{\cal A}_{\bar{F}},\co{1,1}^{(0)}$
by a new one $*$ as follows:
Define three automorphisms $\alpha$, $\beta$, $\gamma$ of 
${\cal A}_{\bar{B}},{\cal A}_{\bar{F}},\co{1,1}^{(0)}$ by
%
%
\begin{equation}
\label{eqn:action}
\alpha(a)\equiv -a,\quad
\beta(f)\equiv -f,\quad
\gamma(s_{i})\equiv (-1)^{i-1}s_{i}\quad(i=1,2).
\end{equation}
Then each of them preserves $\sdag$. 
Define the new involution $*$ on 
${\cal A}_{\bar{B}},{\cal A}_{\bar{F}},\co{1,1}^{(0)}$ by 
\[x^{*}\equiv \alpha(x^{\sdag}),\quad
y^{*}\equiv \beta(y^{\sdag}),\quad
z^{*}\equiv \gamma(z^{\sdag})\quad
(x\in {\cal A}_{\bar{B}},\,y\in{\cal A}_{\bar{F}},\,z\in \co{1,1}^{(0)}).\]
From (\ref{eqn:action}),
we see that (\ref{eqn:abnormalboson}),
(\ref{eqn:abnormalfermion}),
(\ref{eqn:pseudo}) are equivalent to the following equations, respectively:
%
%
\begin{equation}
\label{eqn:abnormalbosontwo}
aa^{*}-a^{*}a=I, \quad aa-aa=a^{*}a^{*}-a^{*}a^{*}=0,
\end{equation}
%
%
\begin{equation}
ff^{*}+f^{*}f=I,\quad ff+ff
=f^{*}f^{*}+f^{*}f^{*}=0,
\label{eqn:abnormalfermiontwo}
\end{equation}
%
%
\begin{equation}
\label{eqn:pseudotwo}
s_{i}^{*}s_{j}=\delta_{ij}I\quad(i,j=1,2),
\quad s_{1}s_{1}^{*}+s_{2}s_{2}^{*}=I.
\end{equation}
New relations (\ref{eqn:abnormalbosontwo}),
(\ref{eqn:abnormalfermiontwo})
and (\ref{eqn:pseudotwo})
are nothing but
CCRs, CARs and the relations of 
canonical generators of the Cuntz algebra $\co{2}$.
Let ${\cal A}_{B}$ and ${\cal A}_{F}$ denote
algebras generated by $\{a,a^{*}\}$ and  $\{f,f^{*}\}$, respectively.
We see that $\alpha,\beta$ and $\gamma$ in (\ref{eqn:action}) are also 
${\bf Z}_{2}$-actions on ${\cal A}_{B}$, ${\cal A}_{F}$ 
and $\co{2}$, respectively
and 
$x^{\sdag}= \alpha(x^{*})$,
$y^{\sdag}= \beta(y^{*})$ and
$z^{\sdag}= \gamma(z^{*})$
for 
$x\in {\cal A}_{B},\,y\in{\cal A}_{F}$ and $z\in \co{2}$.

Furthermore the replacement of involution
is not only the change of appearance but also
compatible to construct representations of 
${\cal A}_{\bar{B}},{\cal A}_{\bar{F}},\co{1,1}^{(0)}$ on
Krein spaces. 
Remark that ${\cal A}_{B}$ and ${\cal A}_{F}$ are same
as ${\cal A}_{\bar{B}}$ and ${\cal A}_{\bar{F}}$ 
as algebras if we take no account of their involutions.

Assume that $({\cal H},\langle\cdot|\cdot\rangle)$ 
is a Hilbert space and $({\cal H},\pi,\eta)$ is a covariant 
representation of the dynamical system
 $({\goth A},{\bf Z}_{2},\varphi)=({\cal A}_{B},{\bf Z}_{2},\alpha),
({\cal A}_{F},{\bf Z}_{2},\beta)$,
$(\co{2},{\bf Z}_{2},\gamma)$,
that is, $\pi\circ \varphi={\rm Ad}\eta\circ \pi$.
Define $(\cdot|\cdot)\equiv \langle\cdot|\eta(\cdot)\rangle$.
Then we see that $\eta$ is a self-adjoint unitary on ${\cal H}$
and $({\cal H},(\cdot|\cdot))$ is a Krein space such that
\[(\pi(x^{\sdag})v|w)=(v|\pi(x)w)\quad(v,w\in {\cal H},\,x\in {\goth A}).\]
Define ${\goth A}_{\pm}\equiv \{x\in {\goth A}: \varphi(x)=\pm x\}$.
Then
%
%
\begin{equation}
\label{eqn:invarianttwo}
\pi({\goth A}_{+}){\cal H}_{\pm}\subset {\cal H}_{\pm},\quad
\pi({\goth A}_{-}){\cal H}_{\pm}\subset {\cal H}_{\mp}
\end{equation}
where ${\cal H}_{\pm}\equiv\{v\in {\cal H}:\eta v=\pm v\}$.
In this way, we see that a covariant representation of 
the involutive algebra $({\goth A},*)$
is closely related to the Krein space representation
of another involutive algebra $({\goth A},\sdag)$.

In Section  \ref{section:second},  we will introduce 
indefinite involutions in order to show
the difference between $\sdag$ and $*$ in the above three examples
and define Krein C$^{*}$-algebras
which are generalizations of C$^{*}$-algebras
including ${\cal A}_{\bar{F}},\co{1,1}^{(0)}$
and the algebra of the Weyl form of abnormal CCRs.
In Sections \ref{section:third}, \ref{section:fourth} and \ref{section:fifth},
we will introduce elementary examples, 
the $\eta$-CCR algebras, the $\eta$-CAR algebras and  the pseudo-Cuntz algebras
as examples of Krein C$^{*}$-algebra.
In Appendix, we will show two models of indefinite-metric quantum field theory.

%
%
\sftt{Krein C$^{*}$-algebra}
\label{section:second}
We introduce indefinite involutions and
Krein C$^{*}$-algebras in this section.
%
%
\subsection{Krein space and indefinite involution on algebra}
\label{subsection:secondone}
A hermitian vector space $({\cal V},(\cdot|\cdot))$ is a {\it Krein space}
if there exists a decomposition ${\cal V}={\cal V}_{+}\oplus {\cal V}_{-}$ 
and $({\cal V}_{\pm},\pm(\cdot|\cdot))$ is a Hilbert space \cite{AI,Bognar,IKL}.
This decomposition is called a {\it fundamental decomposition} of 
$({\cal V},(\cdot|\cdot))$.
By definition,
the new hermitian form $\langle\cdot|\cdot\rangle$ on ${\cal H}$
defined by $\langle v|w\rangle\equiv (v|E_{+}w)-(v|E_{-}w)$ 
for $v,w\in {\cal H}$,
is positive definite 
where $E_{\pm}$ denotes the projection from ${\cal H}$ onto ${\cal H}_{\pm}$.
The operator $U\equiv E_{+}-E_{-}$ is called a {\it fundamental symmetry} 
of $({\cal H},(\cdot|\cdot))$.
For a given Krein space, its fundamental decomposition 
is not unique in general.
Hence we use a {\it Krein triplet} $({\cal H},\langle\cdot|\cdot\rangle,\eta)$,
that is, a Hilbert space $({\cal H},\langle\cdot|\cdot\rangle)$ and
a self-adjoint unitary $\eta$ on ${\cal H}$.
For a Krein triplet $({\cal H},\langle\cdot|\cdot\rangle,\eta)$,
let ${\cal H}_{\pm}\equiv \{v\in {\cal H}:\eta v=\pm v\}$.
Then ${\cal H}={\cal H}_{+}\oplus {\cal H}_{-}$.
Hence $({\cal H},(\cdot|\cdot))$ is a Krein space
with the nondegenerate hermitian form $(\cdot|\cdot)$
defined by $(v|w)\equiv \langle v|\eta w\rangle$ for $v,w\in {\cal H}$.

For an algebra ${\cal A}$ with a unit $I$,
the {\it spectrum} ${\rm sp}_{{\cal A}}(a)$ of $a\in {\cal A}$
is defined by the subset
$\{z\in {\bf C}:\mbox{there exists no inverse of }A-\lambda I
\mbox{ in }{\cal A}\}$ of ${\bf C}$ \cite{Bourbaki}.
We write ${\rm sp}_{{\cal A}}(a)$ as ${\rm sp}(a)$ for simplicity 
of description.
Remark that the definition of the spectrum 
is written without use of any topology.
Important results of spectrum theory do not hold
without use of a norm, especially, the C$^{*}$-condition.
However, there is no assumption of the existence of a norm on
the algebra of field operators in quantum field
theory in general.
Therefore it is necessary to characterize 
a suitable assumption for the involution on such an algebra without topology.
%
%
\begin{defi}
\label{defi:positivity}
Let $({\cal A},*)$ be an involutive algebra with a unit $I$.
\begin{enumerate}
\item
The involution $*$ is positive definite
if $I+x^{*}x$ is invertible for each $x\in {\cal A}$.
\item
The involution $*$ is indefinite if 
there exist $x,y\in {\cal A}$ such that
${\rm sp}(x^{*}x)\cap (0,\infty)\ne\emptyset$ and
${\rm sp}(y^{*}y)\cap (-\infty,0)\ne\emptyset$.
\end{enumerate}
\end{defi}

\noindent
The involution on any C$^{*}$-algebra is positive definite
\cite{Fukamiya,Kaplansky,KelVau}.
If ${\cal A}$ has a unit $I$, then
$I^{*}=I$ for each involution $*$ on ${\cal A}$.
Hence if there exists $x\in {\cal A}$ such that
${\rm sp}(x^{*}x)\cap (-\infty,0)\ne \emptyset$, then $*$ is indefinite.
By definition, if $*$ is positive definite, then
${\rm sp}(x^{*}x)$ is a subset of $[0,\infty)$ for each $x\in {\cal A}$.

%
%
\subsection{Krein C$^{*}$-algebra}
\label{subsection:secondtwo}
We introduce Krein C$^{*}$-algebras in this subsection.
For an involutive algebra $({\cal A},*)$,
we write ${\rm Aut}({\cal A},*)$ 
the set of all involutive automorphisms of $({\cal A},*)$ and 
define ${\rm Aut}_{2}({\cal A},*)\equiv 
\{\alpha \in {\rm Aut}({\cal A},*):\alpha^{2}=id\}$.
For any $\alpha \in{\rm Aut}_{2}({\cal A},*)$,
$\alpha\circ *$ is also an involution on ${\cal A}$.
If $({\cal A},*)$ is a unital C$^{*}$-algebra
and $\alpha\ne id$, then the involution $\sdag$ defined 
by $x^{\sdag}\equiv \alpha(x^{*})$ is indefinite.

We generalize the notion of C$^{*}$-algebra 
according to the definition of Krein space.
%
%
\begin{defi}
\label{defi:krein}
A Banach involutive algebra $({\cal A},\sdag)$ is called 
a Krein C$^{*}$-algebra
if there exists $\alpha\in {\rm Aut}_{2}({\cal A},\sdag)$ such that
\[\|\alpha(x^{\sdag})x\|=\|x\|^{2}\quad \mbox{for all }x\in {\cal A}.\]
In this case,
$\alpha$ is called a fundamental symmetry of $({\cal A},\sdag)$.
\end{defi}

\noindent
By definition, if $({\cal A},\sdag)$ is a Krein C$^{*}$-algebra 
with a fundamental symmetry $\alpha$, 
then $({\cal A},\alpha\circ\sdag)$ is a C$^{*}$-algebra
and $*\equiv \alpha\circ \sdag$ satisfies $*\circ \alpha
=\alpha\circ *$.
We do not know that the uniqueness of fundamental symmetry of
a given Krein C$^{*}$-algebra $({\cal A},\sdag)$.
Hence we define a {\it Krein triplet} $({\cal A},*,\alpha)$
by a C$^{*}$-algebra $({\cal A},*)$
and $\alpha\in {\rm Aut}_{2}({\cal A},*)$. 
This is nothing but a C$^{*}$-dynamical system
$(({\cal A},*),{\bf Z}_{2},\alpha)$.

For two Krein triplets $({\cal A},*,\alpha)$ and $({\cal B},\sdag,\beta)$ 
of C$^{*}$-algebras are {\it isomorphic}
if there exists an involutive isomorphism
$\psi$ from $({\cal A},*)$ and $({\cal B},\sdag)$  such  that
$\psi\circ \alpha=\beta \circ \psi$.
An algebra ${\cal B}$ is a {\it subalgebra}
of Krein triplet $({\cal A},*,\alpha)$ of C$^{*}$-algebra 
if ${\cal B}$ is a C$^{*}$-subalgebra of ${\cal A}$
such that $\alpha({\cal B})={\cal B}$.
If $({\cal A},*,\alpha)$ and $({\cal B},\sdag,\beta)$ are isomorphic,
then $({\cal A},\#,\alpha)$ and $({\cal B},\star,\beta)$ are isomorphic
where $x^{\#}\equiv \alpha(x^{*})$
and $y^{\star}\equiv \beta(y^{\sdag})$ for $x\in {\cal A}$ and $y\in {\cal B}$.

If $({\cal A},*,\alpha)$ is a Krein triplet of C$^{*}$-algebra,
then there is the following natural decomposition as a Banach space:
%
%
\begin{equation}
\label{eqn:decomposition}
{\cal A}={\cal A}_{+}\oplus {\cal A}_{-},\quad
{\cal A}_{\pm}\equiv \{x\in {\cal A}:\alpha(x)=\pm x\}.
\end{equation}
We see that $\{x^{*}:x\in {\cal A}_{\pm}\}\subset {\cal A}_{\pm}$.
Especially, 
${\cal A}_{+}$ is the fixed-point subalgebra ${\cal A}^{\alpha}$ 
of ${\cal A}$ with respect to $\alpha$.
Define $x^{\sdag}\equiv \alpha(x^{*})$ for $x\in {\cal A}$.
Then
\[x^{\sdag}x=\pm x^{*}x \quad(x\in {\cal A}_{\pm}).\]
In this sense,
$\sdag$ on ${\cal A}_{+}$ ({\it resp}. ${\cal A}_{-}$)
is positive definite ({\it resp}. negative definite)
because $x^{*}x\geq 0$ for each $x\in {\cal A}$.
Therefore the decomposition in (\ref{eqn:decomposition})
is regarded as an analogy of a fundamental decomposition of a Krein space.

Because $x^{*}=\alpha(x^{*})$ for each $x\in {\cal A}_{+}$,
the Banach involutive algebra $({\cal A}_{+},\sdag)$ is a C$^{*}$-algebra.
Rewrite $X_{\alpha}\equiv {\cal A}_{-}$.
Then $(X_{\alpha},\langle\cdot|\cdot\rangle)$
is a Hilbert ${\cal A}_{+}$-module
where $\langle a|b\rangle\equiv a^{*}b$ for $a,b\in X_{\alpha}$.
If $x\in {\cal A}$ is $\sdag$-self-adjoint,
then there are unique $x_{\pm}\in {\cal A}_{\pm}$
such that $x_{\pm}^{*}=x_{\pm}$ and $x=x_{+}+\sqrt{-1}x_{-}$.
Hence
\[{\cal A}_{\sdag s.a}=({\cal A}_{+})_{*s.a}
\oplus \sqrt{-1}({\cal A}_{-})_{* s.a}\]
where ${\cal A}_{\sdag s.a}$ denotes the set of all
$\sdag$-self-adjoint elements in ${\cal A}$
and others are defined as the same way.
The part of $\sqrt{-1}({\cal A}_{-})_{* s.a}$
often brings imaginary spectra of $\sdag$-self-adjoint
operators in quantum field theory.

We show the similarity between Krein spaces and Krein C$^{*}$-algebras
as follows:\\

%
%
\noindent

\newcount\secondwid
\secondwid=10

\noindent
\begin{tabular}{lll}
\hline
&
Krein space ${\cal V}$
&Krein C$^{*}$-algebra ${\cal A}$\\
\hline
definition &
\begin{minipage}[c]{4.2cm}
\quad \\
Hilbert space with \\
a unitary
$\eta$, $\eta^{2}=I$
\end{minipage}&
\begin{minipage}[c]{4.5cm}
\quad \\
C$^{*}$-algebra with an\\
automorphism $\alpha$, $\alpha^{2}=id$
\end{minipage}
\\
\hline
structure &
indefinite metric $(\cdot|\cdot)$ &
indefinite involution $\sdag$\\
\hline
\begin{minipage}[c]{1.5cm}
\quad \\
fundamental\\
 symmetry 
\end{minipage}
& $\eta$ & $\alpha$ \\
\hline
\begin{minipage}[c]{1.5cm}
\quad \\
fundamental\\
decomposition
\end{minipage}
&
${\cal V}={\cal V}_{+}\oplus {\cal V}_{-}$ 
&
${\cal A}={\cal A}_{+}\oplus {\cal A}_{-}$
\\
\hline
\begin{minipage}[c]{2.7cm}
\quad \\
positive definite \\
object 
\end{minipage}
&
Hilbert space &
C$^{*}$-algebra\\
\hline
\end{tabular}

%
%
\subsection{Construction of involutive representation of Krein C$^{*}$-algebra}
\label{subsection:secondthree}
We construct an involutive representation of 
an algebra with indefinite involution
from an algebra with positive definite involution.
%
%
\begin{defi}
\label{defi:bounded}
A linear map $T$ on a Krein space $({\cal H},(\cdot|\cdot))$
is bounded if $T$ is bounded with respect to 
the standard Hilbert space of $({\cal H},(\cdot|\cdot))$.
We write ${\cal B}({\cal H})$ the set of all bounded linear operators on ${\cal H}$.
\end{defi}

A triplet $({\cal A},G,\alpha)$ is a {\it C$^{*}$-dynamical system}
if $\alpha$ is a continuous action of the (topological) group $G$ 
on a C$^{*}$-algebra ${\cal A}$.
Especially, if $\alpha\in {\rm Aut}{\cal A}$ satisfies
$\alpha^{2}=id$, then 
we obtain a C$^{*}$-dynamical system $({\cal A},{\bf Z}_{2},\alpha)$ 
associated with $\alpha$. 
%
%
\begin{Thm}
\label{Thm:involutive}
Let $({\cal A},*,\alpha)$ be a Krein triplet of C$^{*}$-algebra 
and define the involution $\sdag$ on ${\cal A}$ by
$x^{\sdag}\equiv \alpha(x^{*})$ for $x\in {\cal A}$.
Assume that ${\cal H}$ is a Hilbert space with 
the positive definite metric $\langle \cdot|\cdot\rangle $ and
$({\cal H},\pi,\eta)$ is a covariant 
representation of the C$^{*}$-dynamical system
$(({\cal A},*),{\bf Z}_{2},\alpha)$.
Define the hermitian form $(\cdot|\cdot)$ on ${\cal H}$ by
%
%
\begin{equation}
\label{eqn:innerproduct}
(v|w)\equiv \langle v|\eta w\rangle \quad(v,w\in {\cal H}).
\end{equation}
Then ${\cal H}$ is a Krein space with respect to
the hermitian form $(\cdot|\cdot)$ and 
$\pi$ is an involutive representation of 
$({\cal A},\sdag)$ on $({\cal H},(\cdot|\cdot))$.
Furthermore the following holds:
\[\pi(x)^{\star}=\eta\pi(x^{*})\eta\quad(x\in{\cal A})\]
where $\star$ denotes
the involution associated with the hermitian form $(\cdot|\cdot)$.
\end{Thm}
%
%
\pr
Since $\pi\circ \alpha={\rm Ad}\eta\circ \pi$,
we see that
\[
\begin{array}{rl}
\langle \pi(x)^{\star}v|w\rangle =&
(\pi(x)^{\star}v|\eta w)\\
 =&(v|\pi(x)\eta w)\\
 =&(v|\eta\pi(\alpha(x))w) \\
=&\langle v|\pi(\alpha(x))w\rangle\\
=&\langle \pi(\alpha(x))^{*}v|w\rangle\\
=&\langle \pi(\alpha(x)^{*})v|w\rangle\\
=&
\langle \pi(\sd{x})v|w\rangle \\
\end{array}
\]
for each $x\in {\cal A}$ and $v,w\in {\cal H}$.
From this, $\pi(x)^{\star}=\pi(\sd{x})$ for each $x\in {\cal A}$.
Hence the former statement holds.
The later is verified by the definition of $\sdag$.
\qedh

\noindent
Under the same assumption in Theorem \ref{Thm:involutive},
we see that 
%
%
\begin{equation}
\label{eqn:invariant}
\pi({\cal A}_{+}){\cal H}_{\pm}\subset {\cal H}_{\pm},\quad
\pi({\cal A}_{-}){\cal H}_{\pm}\subset {\cal H}_{\mp}
\end{equation}
where ${\cal H}_{\pm}\equiv\{v\in {\cal H}:\eta v=\pm v\}$
and ${\cal A}_{\pm}$ is as in (\ref{eqn:decomposition}).
Furthermore
$\pi({\cal A})\cap {\cal B}({\cal H})_{\pm}=\pi({\cal A}_{\pm})$
where ${\cal B}({\cal H})_{+}
\equiv \{a\in {\cal B}({\cal H}):a{\cal H}_{\pm}\subset{\cal H}_{\pm}\}$
and ${\cal B}({\cal H})_{-}
\equiv \{a\in {\cal B}({\cal H}):a{\cal H}_{\pm}\subset{\cal H}_{\mp}\}$.
In this way, the decomposition in (\ref{eqn:decomposition})
is compatible with the fundamental decomposition 
of the Krein space $({\cal H},(\cdot|\cdot))$.
From (\ref{eqn:invariant}) and $*=\sdag$ on ${\cal A}_{+}$,
the Hilbert subspace $({\cal H}_{+},(\cdot|\cdot))$ is an involutive 
representation of the involutive subalgebra ${\cal A}_{+}$
of $({\cal A},\sdag)$.

By Theorem \ref{Thm:involutive}, the following holds.

%
%
\begin{cor}
\label{cor:krein}
Let $({\cal A},*,\alpha)$ be a Krein triplet of C$^{*}$-algebra.
Assume that $\omega$ is a state on $({\cal A},*)$
such that the GNS representation of ${\cal A}$ by $\omega\circ \alpha$ 
is equivalent to that by $\omega$.
Then there exists a self-adjoint unitary $\eta$ on ${\cal H}$
for the GNS representation $({\cal H},\pi)$ by $\omega$
such that $({\cal H},\pi,(\cdot|\cdot))$
is an involutive representation of $({\cal A},\sdag)$
with respect to $(\cdot|\cdot)$ in (\ref{eqn:innerproduct}).
\end{cor}

For a representation $({\cal H},\pi)$ of a C$^{*}$-algebra ${\cal A}$ 
and a C$^{*}$-dynamical system $({\cal A},G,\alpha)$
with a locally compact group $G$,
the {\it regular representation}
$(L_{2}(G,{\cal H}),\tilde{\pi}\cpd \lambda)$
of the crossed product ${\cal A}\cpd G$ 
by $({\cal H},\pi)$ is the representation
which is induced by the following covariant representation
$(L_{2}(G,{\cal H}),\tilde{\pi},\lambda)$ as follows (Section  7.7, \cite{Ped}):
%
%
\begin{equation}
\label{eqn:cov}
(\tilde{\pi}(a)\phi)(g)\equiv \pi (\alpha_{g-1}(a))\phi(g),\quad
(\lambda_{h}\phi)(g)\equiv \phi(h^{-1}g)\quad
\end{equation}
for $a\in {\cal A}$, $g,h\in G$ and $\phi\in L_{2}(G,{\cal H})$
where $L_{2}(G,{\cal H})$ denotes the Hilbert space
of all square integrable  ${\cal H}$-valued functions on $G$
with respect to the Haar measure of $G$.
We apply this for the finite group $G={\bf Z}_{2}$.
Let $({\cal H},\pi)$ be a representation of a C$^{*}$-algebra ${\cal A}$ and
let $\alpha$ be an action of ${\bf Z}_{2}$ on ${\cal A}$.
Assume that there is no unitary $U$ on ${\cal H}$
such that ${\rm Ad}U\circ \pi=\pi\circ \alpha$.
Define the representation $\tilde{\pi}$ of ${\cal A}$ on
the Hilbert space $\tilde{{\cal H}}\equiv {\cal H}\otimes {\bf C}^{2}$ by
%
%
\begin{equation}
\label{eqn:ztwo}
\tilde{\pi}(x)(v\otimes e_{i})\equiv \{\pi(\alpha^{i}(x))v\}\otimes e_{i}
\quad(v\in {\cal H},\,x\in {\cal A},\,i=0,1)
\end{equation}
where $e_{0}$ and $e_{1}$ are standard basis of ${\bf C}^{2}$.
Define the unitary $\eta$ on $\tilde{{\cal H}}$
by 
\[\eta v\otimes e_{i}\equiv v\otimes e_{1-i}\quad(v\in {\cal H},\,i=0,1).\]
Then we can verify that
$\eta\tilde{\pi}(x)\eta^{*}=\tilde{\pi}(\alpha(x))$ for $x\in {\cal A}$.
Hence $(\tilde{{\cal H}},\tilde{\pi},\eta)$
is a covariant representation of the C$^{*}$-dynamical
system $({\cal A},{\bf Z}_{2},\alpha)$.
Define the hermitian form $(\cdot|\cdot)$ on $\tilde{{\cal H}}$ by
$(\cdot|\cdot)\equiv \langle\cdot|\eta(\cdot)\rangle$ and 
$\tilde{{\cal H}}_{\pm}\equiv \{z\in \tilde{{\cal H}}:\eta z=\pm z\}$.
Then
\[\tilde{{\cal H}}_{+}=\{v\otimes (e_{0}+e_{1}):v\in {\cal H}\},\quad
\tilde{{\cal H}}_{-}=\{v\otimes (e_{0}-e_{1}):v\in {\cal H}\}.\]
In this way, we can always construct an involutive representation
which satisfies the assumption in Theorem \ref{Thm:involutive}
from a given involutive representation of a C$^{*}$-algebra
with a ${\bf Z}_{2}$-action.

%
%
\sftt{Elementary examples}
\label{section:third}
We show elementary examples of Theorem \ref{Thm:involutive}.
%
%
\begin{ex}
\label{ex:commutative}
{\rm
Let $C[0,1]$ denote the unital commutative C$^{*}$-algebra
of all complex-valued continuous functions on the interval $[0,1]$ 
with respect to the standard operations.
Define the new involution $\sdag$ on $C[0,1]$ by
\[f^{\sdag}(x)\equiv \overline{f(1-x)}\quad(f\in C[0,1],\,x\in [0,1]).\]
Define $f_{1}\in C[0,1]$ by $f_{1}(x)\equiv 1-2x$.
Then $\{f_{1}^{\sdag}f_{1}\}(x)=-(1-2x)^{2}\leq 0$.
Therefore $\sdag$ is an indefinite involution on $C[0,1]$.
Define $\alpha\in {\rm Aut}C[0,1]$
by $\alpha(f)(x)\equiv f(1-x)$.
Then $f^{\sdag}=\alpha(\bar{f})$.
Define  the representation $\pi$ of $C[0,1]$ 
on the Hilbert space $(L_{2}[0,1],\langle\cdot|\cdot\rangle)$ by
\[\{\pi(f)\phi\}(x)\equiv f(x)\phi(x)\quad 
(f\in C[0,1],\, \phi\in L_{2}[0,1],\,x\in [0,1]).\]
Define the self-adjoint unitary $\eta$ on $L_{2}[0,1]$ by
\[(\eta\phi)(x)\equiv \phi(1-x)\quad(\phi\in L_{2}[0,1],\,x\in [0,1]).\]
Then $(L_{2}[0,1],\pi,\eta)$ is a covariant representation
of the C$^{*}$-dynamical system $(C[0,1],{\bf Z}_{2},\alpha)$.
Hence we obtain the involutive representation
$\pi$ of the Krein C$^{*}$-algebra $(C[0,1],\sdag)$
on the Krein space  $(L_{2}[0,1],(\cdot|\cdot))$
where $(\cdot|\cdot)\equiv \langle \cdot|\eta(\cdot)\rangle$.
}
\end{ex}

%
%
\begin{ex}
\label{ex:matrix}
{\rm (Involutions on $M_{2}({\bf C})$)
The importance of involutions on matrix algebras
are well-known according to Weyl's unitary trick \cite{Weyl,Weyl02,Weyl03}.
We consider relations between such involutions and indefinite 
involutions as follows. 
Let $M_{2}({\bf C})$ denote the unital C$^{*}$-algebra
of all $2\times 2$ matrices with complex entries with respect to the standard operations.
We consider involutions on $M_{2}({\bf C})$
associated with Pauli's spin matrices:
\[
\sigma_{1}=\left(
\begin{array}{cc}
0&1\\
1&0\\
\end{array}
\right),\quad
\sigma_{2}=\left(
\begin{array}{cc}
0&-\sqrt{-1}\\
\sqrt{-1}&0\\
\end{array}
\right)
,\quad
\sigma_{3}=\left(
\begin{array}{cc}
1&0\\
0&-1\\
\end{array}
\right).
\]
Let $\sigma_{0}$ denote the identity matrix for convenience.
Define 
\[\alpha_{i}\equiv {\rm Ad}\sigma_{i}\in {\rm Aut}M_{2}({\bf C}),\quad 
x^{\sdag_{i}}\equiv \alpha_{i}(x^{*})\quad (x\in M_{2}({\bf C}),\,i=0,1,2,3).\]
We rewrite $\ggl(2,{\bf C})=M_{2}({\bf C})$ as a complex Lie algebra with the Lie bracket
of the standard commutator. 
Define the family $\{\gu_{\sdag_{i}}(2):i=0,1,2,3\}$  
of real Lie subalgebras of $\ggl(2,{\bf C})$  by
\[\gu_{\sdag_{i}}(2)\equiv \{X\in \ggl(2,{\bf C}):X^{\sdag_{i}}+X=0\}
\quad(i=0,1,2,3).\]
Then $x^{\sdag_{0}}$ 
is the standard hermitian conjugate of $x\in M_{2}({\bf C})$
and $\gu_{\sdag_{0}}(2)$ is the Lie algebra $\gu(2)$ of $U(2)$.
Let $*\equiv \sdag_{0}$.
Then we can verify that
any two of involutive algebras $(M_{2}({\bf C}),\sdag_{1})$, 
$(M_{2}({\bf C}),\sdag_{2})$ and 
$(M_{2}({\bf C}),\sdag_{3})$ are involutively isomorphic.
Since such isomorphisms can be chosen as preserving the trace, 
$\gu_{\sdag_{1}}(2),\gu_{\sdag_{2}}(2)$ and $\gu_{\sdag_{3}}(2)$ 
are also involutively isomorphic as a Lie algebra.
We explain more concretely as follows:
For
$X=\left(
\begin{array}{cc}
a&b\\
c&d\\
\end{array}
\right)$,
\[\alpha_{1}(X)=
\left(
\begin{array}{cc}
d&c\\
b&a\\
\end{array}
\right),\quad
\alpha_{2}(X)=
\left(
\begin{array}{cc}
d&-c\\
-b&a\\
\end{array}
\right),\quad
\alpha_{3}(X)=
\left(
\begin{array}{cc}
a&-b\\
-c&d\\
\end{array}
\right).
\]
Hence
\[X^{\sdag_{1}}=
\left(
\begin{array}{cc}
\bar{d}&\bar{b}\\
\bar{c}&\bar{a}\\
\end{array}
\right),\quad
X^{\sdag_{2}}=
\left(
\begin{array}{cc}
\bar{d}&-\bar{b}\\
-\bar{c}&\bar{a}\\
\end{array}
\right),\quad
X^{\sdag_{3}}=
\left(
\begin{array}{cc}
\bar{a}&-\bar{c}\\
-\bar{b}&\bar{d}\\
\end{array}
\right).
\]
From these, 
\[
\begin{array}{rl}
\gu_{\sdag_{1}}(2)=&\left\{\left(
\begin{array}{cc}
a&\sqrt{-1}b\\
\sqrt{-1}c&-\bar{a}\\
\end{array}
\right):a\in {\bf C},\,b,c\in {\bf R}\right\},\\
\\
\gu_{\sdag_{2}}(2)=&\ggl(2,{\bf R}),\\
\\
\gu_{\sdag_{3}}(2)=&\{X\in\ggl(2,{\bf C}):I_{1,1}X+X^{*}I_{1,1}
=0\}=\gu(1,1)
\end{array}
\]
where $I_{1,1}=\sigma_{3}$.
We see that  $I_{1,1}^{\sdag_{1}}I_{1,1}=-I$. 
Therefore $\sdag_{1}$ is an indefinite involution on $M_{2}({\bf C})$.
This implies that both $\sdag_{2}$ and $\sdag_{3}$ 
are also indefinite and
neither $\sdag_{1}$,  $\sdag_{2}$ nor $\sdag_{3}$  
is involutively isomorphic to $\sdag_{0}$
because  $\sdag_{0}$ is positive definite.
For any $X\in \gu_{\sdag_{i}}(2)$ $(i=1,2,3)$,
we see that $\exp{X}$ is not unitary.  
However the definition of $\gu_{\sdag_{i}}(2)$
is same as $\gu(2)$ if we misunderstand $\sdag_{i}$ as $*$ for each $i=1,2,3$.

These elementary examples reveal a question about quantum field theory. 
In quantum field theory,
a symmetry is described as a unitary (=anti-hermitian) representation 
of a Lie algebra ${\goth g}$.
Such assumption is formulated as 
\[X^{\sdag}+X=0\]
for a generator $X$ of ${\goth g}$
by using a certain involution $\sdag$. 
However there exists no assumption of the positive definiteness for $\sdag$
because the metric of the representation of 
the theory is unknown until one computes expectation
values of field operators in general.
Therefore we can not know whether $\sdag$ is positive definite or not. 
}
\end{ex}

Examples \ref{ex:commutative} and  \ref{ex:matrix}
show that neither the non-commutativity nor
the infinite-dimensionality of an algebra is essential
for indefinite involution.
The involution $\sdag$ on the algebra generated by
each relation (\ref{eqn:abnormalboson}),
(\ref{eqn:abnormalfermion}) and (\ref{eqn:pseudo}) is indefinite.
%
%
\begin{ex}
\label{ex:etaex}
{\rm
Let $({\cal H},\langle\cdot|\cdot\rangle,\eta)$ be a Krein triplet and
let ${\cal B}({\cal H})$ denote
the C$^{*}$-algebra of all bounded linear operators on 
the Hilbert space $({\cal H},\langle\cdot|\cdot\rangle)$.
For $x\in {\cal B}({\cal H})$, the adjoint 
$x^{\sdag}$ of $x$ with respect to
the hermitian form $(\cdot|\cdot)\equiv \langle\cdot|\eta(\cdot)\rangle$
satisfies
\[x^{\sdag}=\eta x^{*}\eta^{*}.\]
If $\eta\ne I$, then
$\sdag$ is an indefinite involution on ${\cal B}({\cal H})$
because $(\cdot|\cdot)\equiv \langle \cdot|\eta(\cdot)\rangle$
is positive definite if and only if $\sdag$ is positive definite.
}
\end{ex}

%
%
\sftt{$\eta$-CCR and $\eta$-CAR algebra}
\label{section:fourth}
In \cite{IMQ01}, we introduced $\eta$-CCRs and $\eta$-CARs
as families of operators on Krein spaces.
Here we reformulate them without use of representation.
For a Hilbert space $({\cal H},\langle\cdot|\cdot \rangle)$,
let $({\goth A}_{B}({\cal H}),*)$ and $({\goth A}_{F}({\cal H}),*)$
denote the CCR algebra and the CAR algebra over ${\cal H}$,
respectively (Section  5.2.1 in \cite{BraRobi}).
Define $({\cal A}({\cal H}),*)$ the involutive algebra generated by CCRs 
$\{a(f),a^{*}(f):f\in {\cal H}\}$ over ${\cal H}$.
Note that $({\goth A}_{B}({\cal H}),*)$ and 
$({\goth A}_{F}({\cal H}),*)$ are unital C$^{*}$-algebras 
but $({\cal A}({\cal H}),*)$ is not.
Let $\eta$ be a self-adjoint unitary on ${\cal H}$.
%
%
\subsection{$\eta$-CCR algebra}
\label{subsection:fourthone}
Let $({\goth A}_{B}^{(0)}({\cal H},\eta),\sdag)$ denote the involutive algebra
generated by a family $\{W(f):f\in {\cal H}\}$ which satisfies
\[\{W(f)\}^{\sdag}=W(-\eta f),\quad
W(f)W(g)=e^{-\sqrt{-1}{\rm Im}\langle f|g\rangle/2}W(f+g)
\quad(f,g\in {\cal H}).\]
Define the involutive automorphism $\alpha$ on 
$({\goth A}_{B}^{(0)}({\cal H},\eta),\sdag)$ by
$\alpha(W(f))\equiv W(\eta f)$ for $f\in {\cal H}$.
%
%
\begin{lem}
\label{lem:etaccr}
There exists a unique norm $\|\cdot\|$ on ${\goth A}_{B}^{(0)}({\cal H},\eta)$
such that $\|\alpha(x^{\sdag})x\|=\|x\|^{2}$
for each $x\in {\goth A}_{B}^{(0)}({\cal H},\eta)$.
\end{lem}
%
%
\pr
Define the new involution $*$ on ${\goth A}_{B}^{(0)}({\cal H},\eta)$
by $x^{*}\equiv \alpha(x^{\sdag})$ for 
$x\in {\goth A}_{B}^{(0)}({\cal H},\eta)$.
Then we see that $\{W(f):f\in {\cal H}\}$ 
satisfies the canonical relations of the Weyl form of CCRs
with respect to the new involution $*$.
Therefore 
the involutive algebra $({\goth A}_{B}^{(0)}({\cal H},\eta),*)$
is densely embedded into the CCR algebra ${\goth A}_{B}({\cal H})$.
On the other hand, the assumption of the norm in the statement 
is just the C$^{*}$-norm on $({\goth A}_{B}^{(0)}({\cal H},\eta),*)$.
By the uniqueness of the C$^{*}$-norm on 
${\goth A}_{B}({\cal H})$, the statement holds.
\qedh
%
%
\begin{defi}
\label{defi:etaone}
The completion ${\goth A}_{B}({\cal H},\eta)$ of 
${\goth A}_{B}^{(0)}({\cal H},\eta)$ with respect to
the norm in Lemma \ref{lem:etaccr} 
is called the $\eta$-CCR algebra over ${\cal H}$.
\end{defi}
The algebra ${\goth A}_{B}({\cal H},\eta)$
is a Krein C$^{*}$-algebra with a fundamental symmetry $\alpha$.

Let $({\cal A}_{B}({\cal H},\eta),\sdag)$ denote the involutive algebra
generated by a family $\{a(f),a^{\sdag}(f):f\in {\cal H}\}$
which satisfies
\[\left\{
\begin{array}{l}
\{a(f)\}^{\sdag}=a^{\sdag}(f),\\
\\
a(f)a^{\sdag}(g)-a^{\sdag}(g)a(f)=\langle f|\eta g\rangle I,\\
\\
a(f)a(g)-a(g)a(f)=a^{\sdag}(f)a^{\sdag}(g)-a^{\sdag}(g)a^{\sdag}(f)=0\\
\end{array}
\right. \quad(f,g\in {\cal H}).\]
We call $({\cal A}_{B}({\cal H},\eta),\sdag)$ the {\it algebra of 
CCRs} over ${\cal H}$.
Similar algebras are treated in Section  4 of \cite{MMSV}.

%
%
\subsection{$\eta$-CAR algebra}
\label{subsection:fourthtwo}
Let $({\goth A}_{F}^{(0)}({\cal H},\eta),\sdag)$ denote
the involutive algebra generated by a family 
$\{a(f),a^{\sdag}(f):f\in {\cal H}\}$ which satisfies
\[\left\{
\begin{array}{l	}
\{a(f)\}^{\sdag}=a^{\sdag}(f),\\
\\
a(f)a^{\sdag}(g)+a^{\sdag}(g)a(f)=\langle f|\eta g\rangle I,\\
\\
a(f)a(g)+a(g)a(f)=a^{\sdag}(f)a^{\sdag}(g)+a^{\sdag}(g)a^{\sdag}(f)=0\\
\end{array}
\right. \quad(f,g\in {\cal H}).\]
Define the involutive automorphism $\alpha$ on 
$({\goth A}_{F}^{(0)}({\cal H},\eta),\sdag)$ by
$\alpha(a(f))\equiv a(\eta f)$ for $f\in {\cal H}$.
%
%
\begin{lem}
\label{lem:etacar}
There exists a unique norm $\|\cdot\|$ on ${\goth A}_{F}^{(0)}({\cal H},\eta)$
such that $\|\alpha(x^{\sdag})x\|=\|x\|^{2}$
for each $x\in {\goth A}_{F}^{(0)}({\cal H},\eta)$.
\end{lem}
%
%
\pr
In the similarity of the proof of Lemma \ref{lem:etaccr},
the statement holds.
\qedh
%
%
\begin{defi}
\label{defi:etatwo}
The completion ${\goth A}_{F}({\cal H},\eta)$ of 
${\goth A}_{F}^{(0)}({\cal H},\eta)$ with respect to
the norm in Lemma \ref{lem:etacar}
is called the $\eta$-CAR algebra over ${\cal H}$.
\end{defi}

\noindent
The algebra ${\goth A}_{F}({\cal H},\eta)$
is also a Krein C$^{*}$-algebra with a fundamental symmetry $\alpha$.

%
%
\subsection{Representation in Krein space}
\label{subsection:fourththree}
In \cite{IMQ01},
we have already given involutive representations 
($=$ the $\eta$-Fock representations) of 
$({\cal A}_{B}({\cal H},\eta),\sdag)$ and $({\goth A}_{F}({\cal H},\eta),\sdag)$.
%
%
\begin{Thm}
\label{Thm:second}
Let $({\cal H},\langle\cdot|\cdot\rangle,\eta)$ be a Krein triplet and
let ${\cal F}_{+}({\cal H})$ and ${\cal F}_{-}({\cal H})$ denote the
completely symmetric and the completely anti-symmetric Fock space
over the Hilbert space $({\cal H},\langle\cdot|\cdot\rangle)$, respectively.
Let $\Omega$ denote their vacuum vectors as the same symbol.
\begin{enumerate}
\item
There exists a self-adjoint unitary $\Gamma(\eta)$ on ${\cal F}_{+}({\cal H})$
and an involutive representation
$\pi_{B}$ of $({\goth A}_{B}({\cal H},\eta),\sdag)$
on the Krein space $({\cal F}_{+}({\cal H}),(\cdot|\cdot))$ such that
\[(\Omega|\pi_{B}(W(f))\Omega)=e^{-\|f\|^{2}/4}\quad(f\in {\cal H})\]
where $(\cdot|\cdot)$ denotes the hermitian form on ${\cal F}_{+}({\cal H})$
defined by $(v|w)\equiv \langle v|\Gamma(\eta)w\rangle$ 
for $v,w\in {\cal F}_{+}({\cal H})$.
Furthermore, 
$\pi_{B}({\goth A}_{B}({\cal H},\eta))\Omega$ is dense in
${\cal F}_{+}({\cal H})$.
\item
There exists a self-adjoint unitary $\Gamma(\eta)$ on ${\cal F}_{+}({\cal H})$,
a dense subspace ${\cal D}$ of ${\cal F}_{+}({\cal H})$
and an involutive representation
$\pi_{B,0}$ of $({\cal A}_{B}({\cal H},\eta),\sdag)$
on the Krein space $({\cal F}_{+}({\cal H}),(\cdot|\cdot))$ such that
\[\pi_{B,0}(a(f))\Omega=0\quad(\mbox{for all }f\in {\cal H}),\quad
\pi_{B,0}({\cal A}_{B}({\cal H},\eta))\Omega={\cal D}\]
where $(\cdot|\cdot)$ is the hermitian form on ${\cal F}_{+}({\cal H})$
defined by $(v|w)\equiv \langle v|\Gamma(\eta)w\rangle$ 
for $v,w\in {\cal F}_{+}({\cal H})$.
\item
There exists a self-adjoint unitary $\Gamma(\eta)$ on ${\cal F}_{-}({\cal H})$
and an involutive representation $\pi_{F}$ of 
$({\goth A}_{F}({\cal H},\eta),\sdag)$
on the Krein space $({\cal F}_{-}({\cal H}),(\cdot|\cdot))$ such that
\[\pi_{F}(a(f))\Omega=0\quad(\mbox{for all }f\in {\cal H})\]
where $(\cdot|\cdot)$ is the hermitian form on ${\cal F}_{-}({\cal H})$
defined by $(v|w)\equiv \langle v|\Gamma(\eta)w\rangle$ 
for $v,w\in {\cal F}_{-}({\cal H})$.
Furthermore, $\pi_{F}({\goth A}_{F}({\cal H},\eta))\Omega$
is dense in ${\cal F}_{-}({\cal H})$.
\end{enumerate}
Here topologies on ${\cal F}_{+}({\cal H})$ and ${\cal F}_{-}({\cal H})$
are taken as the norm topology induced by the inner product
$\langle\cdot|\cdot\rangle$.
\end{Thm}
%
%
\pr
For (ii) and (iii), see Theorem 1.2 in \cite{IMQ01}.
We show (i).
Let $\Gamma(\eta)$ denote the second quantization of $\eta$ 
and let $\{W(f):f\in {\cal H}\}$ denote the family of Weyl forms of CCRs
on ${\cal F}_{+}({\cal H})$.
Since $\langle\Omega|W(f)\Omega\rangle=e^{-\|f\|^{2}/4}$	
for each $f\in {\cal H}$ in Section  5.2.3 of \cite{BraRobi},
$(\Omega|\pi_{B}(W(f))\Omega)
=\langle\Omega|\Gamma(\eta)\pi_{B}(W(f))\Omega\rangle
=\langle\Omega|\pi_{B}(W(\eta(f))\Gamma(\eta)\Omega\rangle
=\langle\Omega|\pi_{B}(W(\eta(f))\Omega\rangle
=e^{-\|\eta f\|^{2}/4}=e^{-\|f\|^{2}/4}$.
\qedh

%
%
\subsection{Equivalence}
\label{subsection:fourthfour}
When $\eta=I$, the $\eta$-CCRs and the $\eta$-CARs 
coincide with ordinary CCRs and CARs, respectively.
We classify algebras ${\goth A}_{B}({\cal H},\eta)$,
${\cal A}_{B}({\cal H},\eta)$ and ${\goth A}_{F}({\cal H},\eta)$.
%
%
\begin{prop}
\label{prop:equivalence}
For two self-adjoint unitaries $\eta$ and $\eta^{'}$ on ${\cal H}$,
if there exists a unitary $U$ on ${\cal H}$ such that $U\eta U^{*}=\eta^{'}$,
then the following involutive isomorphisms hold:
\[{\goth A}_{B}({\cal H},\eta)\cong {\goth A}_{B}({\cal H},\eta^{'}),
\quad {\cal A}_{B}({\cal H},\eta)\cong {\cal A}_{B}({\cal H},\eta^{'}),
\quad
{\goth A}_{F}({\cal H},\eta)\cong {\goth A}_{F}({\cal H},\eta^{'}).\]
\end{prop}
%
%
\pr
Assume that $\{a(f),a^{\sdag}(f):f\in H\}$ is the set of canonical 
generators of ${\cal A}_{B}({\cal H},\eta)$.
Define $t(f)\equiv a(U^{*}f)$ and $t^{\sdag}(f)\equiv a^{\sdag}(U^{*}f)$ 
for $f\in {\cal H}$.
Then we can verify that $\{t(f),t^{\sdag}(f):f\in {\cal H}\}$ 
satisfy  canonical relations of ${\cal A}_{B}({\cal H},\eta^{'})$.
Hence we obtain the embedding of ${\cal A}_{B}({\cal H},\eta^{'})$ 
into ${\cal A}_{B}({\cal H},\eta)$.
Furthermore this mapping is surjective.
Therefore ${\cal A}_{B}({\cal H},\eta)$ and 
${\cal A}_{B}({\cal H},\eta^{'})$ are involutively isomorphic.
In the same way, we can verify other cases.
\qedh
%
%
\begin{cor}
\label{cor:equ}
In Proposition \ref{prop:equivalence},
equivalences among algebras hold 
if 
\[({\rm ind}_{+}(\eta),{\rm ind}_{-}(\eta))=
({\rm ind}_{+}(\eta^{'}),{\rm ind}_{-}(\eta^{'}))\]
where ${\rm ind}_{\pm}(\eta)\equiv{\rm dim}\{x\in {\cal H}:\eta x=\pm x\}$.
\end{cor}
Since the $\eta$-CCR for ${\rm rank}\eta=2$ 
is unique up to involutive isomorphism (Example 3.3 in \cite{IMQ01}),
the inverse statement of Corollary \ref{cor:equ} does not hold.
%
%
\subsection{Algebra of FP ghosts}
\label{subsection:fourthfive}
We reintroduce the algebra of FP ghosts in \cite{AK05} 
as a Banach involutive algebra.
As for the FP (anti) ghost fields in string theory,  their mode-decomposed 
operators satisfy the abnormal anticommutation relations with the special structure 
owing to the hermiticity of the FP (anti) ghost fields as follows:
%
%
\begin{equation}
\label{eqn:fp}
\left\{
\begin{array}{l}
c_0\bar{c}_0+\bar{c}_0 c_{0} = -I, \quad 
    c_0^{\sdag} = c_0, \ \ \bar{c}_0^{\sdag}=\bar{c}_0, \\
\\
c_m\bar{c}_{n}^{\sdag}+\bar{c}_{n}^{\sdag}c_m
 = c_m^{\sdag}\bar{c}_n+\bar{c}_n c_m^{\sdag} =-\delta_{m,n}I 
    \quad (m, \, n= 1,\,2\,\ldots)                
\end{array}
\right.
\end{equation}
and other anticommutation relations vanish.
Define ${\cal F}{\cal P}_{0}$ 
the involutive algebra generated by $\{c_{n},\,\bar{c}_{n}:n\geq 0\}$.
Define the self-adjoint unitary $\eta$ on 
the Hilbert space ${\cal H}\equiv l_{2}({\bf Z}_{\geq 0})$ by
%
%
\begin{equation}
\label{eqn:fpeta}
\eta e_{2n}\equiv -e_{2n+1},\quad
\eta e_{2n+1}\equiv -e_{2n}\quad(n\geq 0)
\end{equation}
where ${\bf Z}_{\geq 0}\equiv \{n\in {\bf Z}:n\geq 0\}$.
%
%
\begin{lem}
\label{lem:extend}
For $\eta$ in (\ref{eqn:fpeta}),
there exists an involutive embedding 
of ${\cal F}{\cal P}_{0}$ into the $\eta$-CAR
algebra ${\goth A}_{F}({\cal H},\eta)$.
\end{lem}
%
%
\pr
Let $\{a(f):f\in {\cal H}\}$ denote the canonical 
generators of ${\goth A}_{F}({\cal H},\eta)$ and define
$a_{n}\equiv a(e_{n})$ for $n\geq 0$
where $\{e_{n}\}_{n\geq 0}$ denotes the standard basis of ${\cal H}$.
Define the map $\varphi$ from 
${\cal F}{\cal P}_{0}$ into ${\goth A}_{F}({\cal H},\eta)$ by
\[
\begin{array}{l}
\varphi(c_{0})\equiv 2^{-1/2}(a_{0}+a_{0}^{\sdag}),\quad
\varphi(\bar{c})\equiv 2^{-1/2}(a_{1}+a_{1}^{\sdag}),\\
\\
\varphi(c_{n})\equiv a_{2n},\quad
\varphi(\bar{c}_{n})\equiv a_{2n+1}\quad(n\geq 1).
\end{array}
\]
Then we can verify that $\varphi$ is an involutive embedding.
\qedh

\noindent
The completion ${\cal F}{\cal P}$ 
of ${\cal F}{\cal P}_{0}$ in ${\goth A}_{F}({\cal H},\eta)$ 
is the Krein C$^{*}$-algebra of FP ghosts.

%
%
\sftt{$\eta$-Cuntz algebra}
\label{section:fifth}
In this section, we introduce $\eta$-Cuntz algebras
as Krein C$^{*}$-algebras.
%
%
\subsection{Definition and equivalence}
\label{subsection:fifthone}
Let $\eta=(\eta_{ij})_{i,j=1}^{N}$ be a self-adjoint unitary in $U(N)$. 
Let $(\co{\eta}^{(0)},\sdag)$ denote the unital algebra $\co{\eta}^{(0)}$
with an involution $\sdag$
generated by $s_{1},\ldots,s_{N}$ which satisfies
%
%
\begin{equation}
\label{eqn:etacuntz}
s^{\sdag}_{i}s_{j}=\eta_{ij}I\quad(i,j=1,\ldots,N),\quad
\sum_{i,j=1}^{N}\eta_{ij} s_{i}s^{\sdag}_{j}=I.
\end{equation}
Define the involutive automorphism $\alpha_{\eta}$ on 
$(\co{\eta}^{(0)},\sdag)$ by 
%
%
\begin{equation}
\label{eqn:etaauto}
\alpha_{\eta}(s_{i})\equiv \sum_{j=1}^{N}\eta_{ji}s_{j}\quad(i\edot).
\end{equation}
%
%
\begin{lem}
\label{lem:norm}
The algebra $\co{\eta}^{(0)}$ has a unique norm $\|\cdot\|$ such that
$\|\alpha_{\eta}(x^{\sdag})x\|=\|x\|^{2}$ for each $x\in \co{\eta}^{(0)}$.
\end{lem}
%
%
\pr
Define the new involution $*$ on $\co{\eta}^{(0)}$ by
$x^{*}\equiv \alpha_{\eta}(x^{\sdag})$ for $x\in \co{\eta}^{(0)}$.
Then we see that $(\co{\eta}^{(0)},*)$ is a dense 
involutive subalgebra of the Cuntz algebra $(\con,*)$ \cite{C}.
Since the norm $\|\cdot\|$ satisfies 
the C$^{*}$-condition with respect to $*$,
the norm $\|\cdot\|$ is a unique C$^{*}$-norm on the involutive algebra
$(\co{\eta}^{(0)},*)$.
\qedh

\noindent
The existence of $\co{\eta}^{(0)}$ for any $\eta$
is also verified in the proof of Lemma \ref{lem:norm}.

%
%
\begin{defi}
\label{defi:eta}
\begin{enumerate}
\item
The completion $\co{\eta}$ of $\co{\eta}^{(0)}$ with respect to
the norm in Lemma \ref{lem:norm} is called the $\eta$-Cuntz algebra.
\item
For nonnegative integers $d,d^{'}$ with $d+d^{'}\geq 2$,
the Banach involutive algebra $(\co{d,d^{'}},\sdag)$
is called the pseudo-Cuntz algebra
if $\co{d,d^{'}}$ is the $\eta$-Cuntz algebra for 
$\eta\equiv (\eta_{ij})_{i,j=1}^{d+d^{'}}\in U(d+d^{'})$ defined by
%
%
\begin{equation}
\label{eqn:eta}
(\eta_{ij})_{i,j=1}^{d+d^{'}}\equiv 
\left(\begin{array}{cc}
I_{d}&0\\
0& -I_{d^{'}}\\
\end{array}
\right).
\end{equation}
\end{enumerate}
\end{defi}

\noindent
Especially, $\co{d,0}\cong \co{d}$.
In \cite{AK05}, we generalized the Cuntz algebra
to the \pca\ without topology which is an involutive algebra 
of operators on an indefinite-metric space.
In Definition \ref{defi:eta} (ii), 
the \pca\ is redefined as a Banach involutive algebra
without use of any representation.

The equivalence among $\{\co{\eta}:\eta\in U(N),\,\eta^{*}=\eta\}$ 
is shown as follows:
%
%
\begin{lem}
\label{lem:lambda}
If $\eta=\Lambda \eta^{'}\Lambda^{*}$ for a unitary $\Lambda\in U(N)$,
then $\co{\eta}\cong \co{\eta^{'}}$.
\end{lem}
%
%
\pr
Let $s_{1},\ldots,s_{N}$ denote the canonical generators of $\co{\eta}$.
Define elements 
$t_{i}\equiv \sum_{j=1}^{N}\Lambda^{*}_{ji}s_{j}$
in $\co{\eta^{'}}$  for $i=1,\ldots,N$.
Then we see that $t_{1},\ldots,t_{N}$ satisfy 
the canonical relations of $\co{\eta^{'}}$.
Therefore $\co{\eta^{'}}$ is embedded into $\co{\eta}$.
Furthermore this embedding is surjective.
Hence the statement holds.
\qedh

%
%
\begin{cor}
\label{cor:self}
For any self-adjoint element $\eta\in U(N)$,
the algebra $\co{\eta}$ is isomorphic to a certain pseudo-Cuntz algebra.
\end{cor}
Let $({\cal A},\sdag)$ denote the involutive algebra generated by
$s_{1}$ and $s_{2}$ which satisfy the following:
\[s_{1}^{\sdag}s_{2}=I,\quad s_{1}^{\sdag}s_{1}=s_{2}^{\sdag}s_{2}=0,\quad
s_{1}s_{2}^{\sdag}+s_{2}s_{1}^{\sdag}=I.\]
Then the involution $\sdag$ is indefinite.
The algebra $({\cal A},\sdag)$ is densely, involutively embedded 
into $(\co{1,1},\sdag)$.
%
%
\begin{ex}
\label{ex:cuntz}
{\rm
Two algebras $\co{2}$ and $\co{1,1}$ are not involutive isomorphic
because there is no nondegenerate  involutive representation 
of $\co{1,1}$ on Hilbert space.
In the same reason,
$\co{2}$ and $\co{0,2}$ are not involutive isomorphic.
We let the following open problems:
Whether are $\co{2,1}$ and $\co{1,2}$ equivalent or not?
Whether are $\co{1,1}$ and $\co{0,2}$ equivalent or not?
Classify all pseudo-Cuntz algebras.
}
\end{ex}

For $\eta$ in (\ref{eqn:eta}),
$\rho$ is the {\it canonical endomorphism} of $\co{d,d^{'}}$
if $\rho$ is the map on $\co{d,d^{'}}$ defined by 
\[\rho(x)\equiv \sum_{i=1}^{d+d^{'}}\eta_{ii} s_{i}xs_{i}^{\sdag}.\]
We see that $\rho(x)\rho(y)=\rho(xy)$ for $x,y\in \co{d,d^{'}}$.

For $\eta$ in (\ref{eqn:eta}),
define $U(d,d^{'})\equiv\{g\in GL_{d+d^{'}}({\bf C}):g\eta g^{*}=\eta\}$.
Then $U(d,d^{'})$ is a group such that
$g^{*}\in U(d,d^{'})$ for each $g\in U(d,d^{'})$
where $*$ denotes the hermite conjugate on $M_{d+d^{'}}({\bf C})$.
We see that 
$U(d,d^{'})=\{g\in GL_{d+d^{'}}({\bf C}):g^{*}\eta g=\eta\}$.
For $g\in U(d,d^{'})$, define the involutive
automorphism $\alpha_{g}$ of $\co{d,d^{'}}$ by
\[\alpha_{g}(s_{i})\equiv \sum_{j=1}^{d+d^{'}}g_{ji}s_{j}\quad(i=1,\ldots,d+d^{'}).\]
Then $\alpha$ is an involutive action of $U(d,d^{'})$ on $(\co{d,d^{'}},\sdag)$.
Especially, the $U(1)$-gauge action on $\co{d,d^{'}}$ is also an involutive
action.

%
%
\subsection{Involutive representation of $\co{d,d^{'}}$ on Krein space}
\label{subsection:fifthtwo}
According to Theorem \ref{Thm:involutive},
we consider involutive representations of pseudo-Cuntz algebras.
%
%
\begin{cor}
\label{cor:assume}
Let $\eta$ be as in (\ref{eqn:eta}).
Assume that $({\cal H},\pi,U)$ 
is the covariant representation of the C$^{*}$-dynamical system
$(\co{d+d^{'}},{\bf Z}_{2},\alpha_{\eta})$
for $\alpha_{\eta}$ in (\ref{eqn:etaauto}) and $N=d+d^{'}$.
Define the hermitian form $(\cdot|\cdot)$ on ${\cal H}$ by
\[(v|w)\equiv \langle v |Uw\rangle \quad(v,w\in {\cal H}).\]
Then $\pi$ is an involutive representation of $(\co{d,d^{'}},\sdag)$
on the Krein space $({\cal H},(\cdot|\cdot))$.
\end{cor}

We show examples.
%
%
\begin{ex}
\label{ex:cuntztwo}
{\rm
Let $({\cal H},\pi)$ be a representation of $\co{2}$ 
with a cyclic vector $\Omega\in {\cal H}$
which satisfies 
\[\pi(s_{1}s_{2})\Omega=\Omega.\]
Such $({\cal H},\pi)$ is $P(12)$ in \cite{IWF01} and 
it is unique up to unitary equivalence and irreducible.
Define $\alpha\in {\rm Aut}\co{2}$ and the new involution
$\sdag$ on $\co{2}$ by 
\[\alpha(s_{1})\equiv s_{1},\quad \alpha(s_{2})\equiv -s_{2},\quad
x^{\sdag}\equiv \alpha(x^{*})\quad(x\in \co{2}).\]
Then $\co{2}$ becomes  $\co{1,1}$ by replacing the involution $*$
with $\sdag$.
We construct an involutive representation
of $(\co{1,1},\sdag)$ from $({\cal H},\pi)$ as follows.
Since $({\cal H},\pi\circ \alpha)$ is not equivalent to
$({\cal H},\pi)$ ,
$({\cal H},\pi)$ itself is not a covariant representation 
of $(\co{2},{\bf Z}_{2},\alpha)$.
From $({\cal H},\pi)$,
we obtain the representation
$(\tilde{{\cal H}},\tilde{\pi})$ of $\co{2}$ in (\ref{eqn:ztwo}).
Define $V_{1}\equiv \tilde{\pi}(\co{2})\Omega\otimes {\bf C}e_{0}$
and $V_{2}\equiv \tilde{\pi}(\co{2})\Omega\otimes {\bf C}e_{1}$.
We see that both $\tilde{\pi}|_{V_{1}}$
and $\tilde{\pi}|_{V_{2}}$
are irreducible representations of $\co{2}$.
On the other hand,
${\cal H}_{+}\equiv\{v\otimes (e_{0}+e_{1}):v\in {\cal H}\}$
and ${\cal H}_{-}\equiv\{v\otimes (e_{0}-e_{1}):v\in {\cal H}\}$.
Hence
${\cal H}_{\pm}$ is not invariant under the action of $\tilde{\pi}(\co{1,1})$.
}
\end{ex}

%
%
\begin{ex}
\label{ex:cuntzthree}
{\rm
Let $d\geq 1$ and $\eta_{ij}\equiv (-1)^{i-1}\delta_{ij}$ for $i,j=1,2,\ldots,2d$.
Define the $*$-automorphism of $\co{2d}$ by
%
%
\begin{equation}
\label{eqn:alpha}
\alpha(s_{i})\equiv \eta_{ii} s_{i}\quad(i=1,2,3,\ldots, 2d).
\end{equation}
Then $\alpha^{2}=id$. Define the new involution
$\sdag$ on $\co{2d}$ by $x^{\sdag}=\alpha(x^{*})$ for $x\in \co{2d}$.
Then we see that 
\[s_{i}^{\sdag}s_{j}=\eta_{ij}I\quad(i,j=1,2,3,\ldots,2d),\quad
\sum_{i,j=1}^{2d}\eta_{ij}s_{i}s_{j}^{\sdag}=I.\]
Hence they generate $\co{d,d}$.
In this way,
the replacement of involution on $\co{2d}$ makes $\co{2d}$ to $\co{d,d}$.

Next, we construct an involutive representation of $\co{d,d}$ as follows.
Let $(({\cal H},\langle\cdot|\cdot\rangle),\pi)$ be a $*$-representation of $\co{2d}$
with a cyclic unit vector $\Omega\in {\cal H}$
which satisfies 
\[\pi(s_{1})\Omega=\Omega.\]
By this assumption, we see that $({\cal H},\pi)$ is irreducible.
Let $\{1,\ldots,2d\}^{*}\equiv \bigcup_{k\geq 0}\{1,\ldots,2d\}^{k}$
where $\{1,\ldots,2d\}^{0}\equiv \emptyset$.
Define the function $\chi$ on $\{1,\ldots,2d\}^{*}$ by
$\chi(\emptyset)\equiv 1$,
\[\chi(J)\equiv (-1)^{n_{2}(J)}\]
where $n_{2}(J)\equiv \#\{i\in\{1,\ldots,k\}:j_{i}\mbox{ is even}\}$
for $J=(j_{1},\ldots,j_{k})$.
Define the subset $\Lambda\equiv \{(1),J\cup (i):J\in\{1,2,\ldots,2d\}^{*},\,
2\leq i\leq 2d\}$ of multiindices and
let subsets $\Lambda_{\pm}$ of $\Lambda$ by
$\Lambda_{\pm}\equiv \{J\in\Lambda:\chi(J)=\pm 1\}$.
Define two subspaces ${\cal H}_{\pm,0}$ of ${\cal H}$ by
\[{\cal H}_{\pm,0}\equiv {\rm Lin}\langle\{\pi(s_{J})\Omega:J\in \Lambda_{\pm}\}
\rangle\]
and let ${\cal H}_{\pm}$ denote their completions.
Then $\{\pi(s_{J})\Omega: J\in\Lambda\}$ is a complete
orthonormal basis of ${\cal H}$
and ${\cal H}={\cal H}_{+}\oplus {\cal H}_{-}$.
Define the unitary $\eta$ on ${\cal H}$ by
\[\eta\pi(s_{J})\Omega\equiv \chi(J)\pi(s_{J})\Omega\quad(J\in\Lambda).\]
For $\alpha$ in (\ref{eqn:alpha}),
we can verify that 
$\pi\circ \alpha={\rm Ad}\eta\circ \pi$.
Let $e_{1}\equiv \Omega$ and $e_{i}\equiv \pi(s_{i})\Omega$ 
for $i=2,3,4,\ldots,2d$
and define $\{e_{n}\in {\cal H}:n\in {\bf N}\}$ recursively by
\[e_{4d(n-1)+i}\equiv \pi(s_i)e_{2n-1}, \quad
e_{4dn + 1 -i}\equiv \pi(s_{i})e_{2n}\quad(i=1,2,\ldots,2d,,n\in {\bf N}).\]
Then $\{e_{2n-1}:n\in {\bf N}\}\subset{\cal H}_{+}$ and
$\{e_{2n}:n\in {\bf N}\}\subset{\cal H}_{-}$.
Define the new hermitian form $(\cdot|\cdot)$ on ${\cal H}$ by
\[(v|w)\equiv \langle v|\eta w\rangle\quad(v,w\in {\cal H}).\]
Then $(e_{n}|e_{m})=\eta_{nm}$ for $n,m\in {\bf N}$.
We see that $({\cal H},(\cdot|\cdot))$ is a Krein space
with a fundamental decomposition
${\cal H}_{+}\oplus {\cal H}_{-}$
and $(({\cal H},(\cdot|\cdot)),\pi)$ is an involutive representation of 
the pseudo-Cuntz algebra $(\co{d,d},\sdag)$.
This is just the example in Section  3 of \cite{AK05}.
}
\end{ex}


\ssfr{Acknowledgement}

The author would like to express his sincere thanks to 
Noboru Nakanishi, Mitsuo Abe and Takeshi Nozawa
for their advices on this subject.

\appendix

%
%
\sftt{Pauli's modified Schr\"{o}dinger representation of 
abnormal canonical commutation relations}
\label{section:appone}
We review the Pauli's example in Section  3 of \cite{Pauli}
as a covariant representation of the involutive algebra ${\cal A}$ generated
by abnormal canonical commutation relations 
by modifying the Schr\"{o}dinger representation.
He strictly distinguished ``hermitian operator" and ``self-adjoint operator",
and $``*"$ and $``\sdag"$ as his notation and terminology.
Remark that we change the notation and terminology from originals.
Let $({\cal A},\sdag)$ denote the involutive algebra
generated by $a$ and $a^{\sdag}$ which satisfy
(\ref{eqn:abnormalboson}).
We see that ${\cal A}={\cal A}_{B}({\bf C},-I)$ 
in Section  \ref{subsection:fourthone}.
Define the involutive automorphism $\alpha$ of $({\cal A},\sdag)$ by
$\alpha(a)\equiv -a$ and $\alpha(a^{\sdag})\equiv -a^{\sdag}$.
Define the new involution $*$ on ${\cal A}$ by
$x^{*}\equiv \alpha(x^{\sdag})$ for $x\in {\cal A}$.
Then we see that $aa^{*}-a^{*}a=I$.
Define the covariant representation 
$({\cal H},\pi,R)$ of the involutive dynamical system
$(({\cal A},*),{\bf Z}_{2},\alpha)$ by
\[{\cal H}\equiv L_{2}({\bf R}),\quad
\pi(a)\equiv 2^{-1/2}(\hat{p}-\sqrt{-1}\hat{q}),\]
\[R:{\cal H}\to {\cal H};\quad
(Rf)(q)\equiv f(-q)\quad(q\in {\bf R})\]
where $\hat{p}\equiv -\sqrt{-1}d/dq$
and let $\hat{q}$ denote the position operator on ${\cal H}$.
Since operators $\hat{p}$ and $\hat{q}$ are 
extended to self-adjoint operators on ${\cal H}$,
the adjoint $\pi(a)^{*}$ of $\pi(a)$ 
with respect to the standard inner product on $L_{2}({\bf R})$ is given by
$\pi(a)^{*}= 2^{-1/2}(\hat{p}+\sqrt{-1}\hat{q})$.
We see that $\pi\circ \alpha={\rm Ad}R\circ \pi$.
Hence the algebra $({\cal A},\sdag)$ is involutively represented on
the Krein space $({\cal H},(\cdot|\cdot))$ where
the indefinite metric $(\cdot|\cdot)$ on ${\cal H}$ is defined by
\[(f|g)\equiv 
\int _{{\bf R}}\overline{f(q)}\,\{Rg\}(q)\,dq=
\int _{{\bf R}}\overline{f(q)}\,g(-q)\,dq\quad
(f,g\in {\cal H}).\]
Then the adjoint $\pi(a)^{\sdag}$ of $\pi(a)$
with respect to $(\cdot|\cdot)$ is given by
\[\pi(a)^{\sdag}= 2^{-1/2}(-\hat{p}-\sqrt{-1}\hat{q})=-\pi(a)^{*}.\]
Furthermore $\hat{p}^{\sdag}=-\hat{p}$ and $\hat{q}^{\sdag}=-\hat{q}$.
For ${\cal H}_{\pm}\equiv \{v\in {\cal H}:R v=\pm v\}$,
we see that ${\cal H}_{+}$ ({\it resp}. ${\cal H}_{-}$) is the 
space of all even ({\it resp}. odd) functions in ${\cal H}$.
In consequence, the replacement of involution
on the Schr\"{o}dinger representation gives an involutive representation
of abnormal canonical commutation relations on the Krein space.

%
%
\sftt{A model with indefinite metric}
\label{section:apptwo}
It is known that there are various simple models associated with
indefinite-metric quantum field theory \cite{Araki,Ito,Nakanishi}.
Araki treated a system of a boson and an abnormal boson \cite{Araki}.
Nakanishi treated a system of a fermion, an abnormal fermion and 
a countably infinite family of bosons in the Lee model 
(Section 12, \cite{Nakanishi}. See also \cite{Fuda}).
They computed eigenvalues of Hamiltonians which are apparently self-adjoint,
and showed that there exist complex eigenvalues which are not real
under certain conditions.
Such Hamiltonian is treated 
as ``pseudo-Hermitian operator" in \cite{Solombrino}.

In order to  effectively explain how the indefinite-metric space appears
and why non-real valued eigenvalues are derived,
we introduce a simpler model 
as a {\it virtual system} of transformation among 
a particle $A$ to a particle $B$
\[A\leftrightarrows B.\]
Assume that $A$ is a boson ($=b$) or a fermion ($=f$),
and $B$ is an abnormal boson ($=\bar{b}$) or an abnormal fermion ($=\bar{f}$).
Hence there are four combinations of the choice of particles $A$ and $B$
as follows:
%
%
\begin{equation}
\label{eqn:combination}
(A,B)=(b,\bar{b}),(b,\bar{f}),(f,\bar{b}),(f,\bar{f}).
\end{equation}
For every combination, we define the (common) Hamiltonian $H$ by
%
%
\begin{equation}
\label{eqn:hamiltonian}
H\equiv m_{A}a_{A}^{\sdag}a_{A}-m_{B}a_{B}^{\sdag}a_{B}
+ga_{A}^{\sdag}a_{B}+\bar{g}a_{B}^{\sdag}a_{A}
\end{equation}
where $m_{X},a^{\sdag}_{X},a_{X}$ denote the mass,
the creation and the annihilation operator of $X=A,B$
and $g \in {\bf C}$ is their coupling constant.
We assume that 
$a_{A}a_{B}=a_{B}a_{A}$ and $a_{A}^{\sdag}a_{B}=a_{B}a_{A}^{\sdag}$.
Let $\Omega$ denote the common vacuum vector
such that $a_{A}\Omega=a_{B}\Omega=0$.
Then we see that the hermitian form $(\cdot|\cdot)$ on the state space such that
$H$ is self-adjoint and $\Omega$ is normalized,  satisfies the following:
%
%
\begin{equation}
\label{eqn:indefinite}
(\Omega|\Omega)=(a^{\sdag}_{A}\Omega|a^{\sdag}_{A}\Omega)=1,
\quad
(a^{\sdag}_{B}\Omega|a^{\sdag}_{B}\Omega)=-1
\end{equation}
and $\Omega,a^{\sdag}_{A}\Omega,
a^{\sdag}_{B}\Omega$ are mutually orthogonal.
%
%
\begin{prop}
\label{prop:unifying}
Let $V\equiv {\rm Lin}\langle\{a_{A}^{\sdag}\Omega,a_{B}^{\sdag}\Omega\}\rangle$
and let $H$ be as in (\ref{eqn:hamiltonian}).
Then $HV\subset V$ and the eigenvalue of $H$ on $V$ is the following:
\begin{enumerate}
\item
If $|m_{A}-m_{B}|>2|g|$,
then $H$ has two different real eigenvalues.
The norm of one of these two eigenvectors is negative.
\item
If $|m_{A}-m_{B}|=2|g|$,
then $H$ has one real eigenvalue.
The associated eigenvector is neutral, that is, zero-norm.
\item
If $|m_{A}-m_{B}|<2|g|$,
then $H$ has two different eigenvalues which are not real.
\end{enumerate}
\end{prop}
%
%
\pr
For each case in (\ref{eqn:combination}),
the first statement is easily verified and
the matrix representation of $H$ with respect to
$a^{\sdag}_{A}\Omega$ and $a^{\sdag}_{B}\Omega$ is 
%
%
\begin{equation}
\label{eqn:restriction}
H|_{V}=\left(
\begin{array}{cc}
m_{A}& -g\\
\bar{g}&m_{B}\\
\end{array}
\right).
\end{equation}
By the discriminant of the eigenequation of this matrix, 
the statement (i),(ii) and (iii) about the properties of eigenvectors hold.
The norm of the eigenvector with respect to the indefinite metric
$(\cdot|\cdot)$ is verified by direct computation.
\qedh

\noindent
The operator $H$ in (\ref{eqn:hamiltonian})
is hermite on the state space with
the indefinite metric $(\cdot|\cdot)$
which satisfies (\ref{eqn:indefinite}).
However the matrix in (\ref{eqn:restriction}) is not hermite
where (\ref{eqn:restriction}) is obtained from
the comparison of the left side and the right side of the eigenequation.
The result of Proposition \ref{prop:unifying}
is similar to that of the Lee model in \cite{Nakanishi}.
For our model, it is not necessary to use the topology of 
indefinite-metric space and operators on it because ${\rm dim}V<\infty$.

We consider the above model as a representation theory of involutive algebra.
Let $({\cal A}_{X},\sdag)$ denote the involutive algebra generated by $a_{X}$
and let ${\cal H}_{X}\equiv {\cal F}_{+}({\bf C},\eta)$ or
${\cal F}_{-}({\bf C},\eta)$ for $X=A,B$, respectively.
If $A$ is a boson, then  ${\cal H}_{A}\cong \ltn$ and $\eta=I$.
If $A$ is a fermion, then  ${\cal H}_{A}\cong {\bf C}^{2}$ and $\eta=I$.
The total algebra is ${\cal A}_{A}\otimes {\cal A}_{B}$
and the representation space is ${\cal H}_{A}\otimes {\cal H}_{B}$.
The Hamiltonian $H$ belongs to ${\cal A}_{A}\otimes {\cal A}_{B}$.
Define the involutive automorphism $\alpha$ of ${\cal A}_{A}\otimes {\cal A}_{B}$ by
$\alpha(a_{A})\equiv a_{A}, \alpha(a_{A}^{\sdag})\equiv a_{A}^{\sdag},
\alpha(a_{B})\equiv -a_{B}$ and $\alpha(a_{B}^{\sdag})\equiv -a_{B}^{\sdag}$
where we identify $a_{A}$ and $a_{B}$ with $a_{A}\otimes I$
and $I\otimes a_{B}$, respectively.
For the new involution $x^{*}\equiv \alpha(x^{\sdag})$,
we obtain that 
\[H=m_{A}a_{A}^{*}a_{A}+m_{B}a_{B}^{*}a_{B}
+ga_{A}^{*}a_{B}-\bar{g}a_{B}^{*}a_{A}.\]
Then we see that $H$ is not self-adjoint with respect to $*$ on
the positive definite-metric space.
Since the eigenvalue of $H$ is independent in the choice of
hermitian form, the result in (\ref{prop:unifying}) holds.


%

\begin{thebibliography}{99}
%
%
\bibitem{AK01}Abe, M.,  and Kawamura, K.:
Recursive fermion system in Cuntz algebra.\ I 
---Embeddings of fermion algebra into Cuntz algebra---,
Comm. Math. Phys. {\bf 228} (2002) 85-101.
%
\bibitem{AK05}Abe, M.  and Kawamura, K.:
Pseudo-Cuntz algebra and recursive FP ghost system in string theory,
Int. J. Mod. Phys. {\bf A18}, No. 4 (2003) 607-625.
%
\bibitem{Araki}Araki, H.:
On a pathology in indefinite metric inner product space, 
Comm. Math. Phys. {\bf 85} (1982), no. 1, 121--128.
%
\bibitem{AI}Azizov, T. Ya.  and Iokhvidov, I. S.: 
{\it  Linear operators in spaces with an indefinite metric}, 
John Wiley \& Sons, Chichester, 1989.
%
\bibitem{Bognar}Bogn\'ar, J.: 
{\it Indefinite inner product spaces},
Springer-Verlag, New York, 1974.
%
\bibitem{Bourbaki}Bourbaki, N.: 
{\it \'{E}l\'{e}ments de math\'{e}matique, th\'{e}ories spectrales},
Hermann Paris, 1967.
%
\bibitem{BraRobi}Bratteli, O.  and Robinson, D. W.: 
{\it Operator algebras and quantum statistical mechanics II},
Springer New York (1981).
%
\bibitem{C}Cuntz, J.: 
Simple $C^*$-algebras generated by isometries,
Comm. Math. Phys. {\bf 57} (1977) 173-185.
%
\bibitem{Dirac}Dirac, P. M.: 
The physical interpretation of quantum mechanics,
Proc. Roy.  Soc. London Ser. A {\bf 180} (1942) 1--40.  
%
\bibitem{Fuda}Fuda,  M. G.: 
Poincare invariant Lee model,
Phys. Rev. D {\bf 41} 2 (1990) 534--549.
%
\bibitem{Fukamiya}Fukamiya, M.: 
On a theorem of Gel'fand and Naimark and the B$^{*}$-algebra,
 Kumamoto J. Sci. {\bf 1} (1952), 17--22.
%
%
\bibitem{IKL}Iohvidov, I. S.,  Krein, M. G.  and Langer, H.: 
{\it Introduction to the spectral theory of operators in spaces with an
indefinite metric},
 Mathematical Research, 9. Akademie-Verlag, Berlin, 1982.
%
\bibitem{Ito}Ito, K. R.: 
Canonical linear transformation on Fock space with an indefinite metric,
Publ. RIMS Kyoto Univ. {\bf 14} (1978) 503--556.
%
\bibitem{Kaplansky}Kaplansky, I.: 
Cited by J. A. Schultz in a review of [Fuk 1],
 Math. Rev. {\bf 14} (1953), 884.
%
\bibitem{IWF01}Kawamura, K.: 
Extensions of representations of the CAR algebra to
the Cuntz algebra ${\cal O}_2$ ---the Fock and the infinite wedge---,
J. Math. Phys. {\bf 46}, no. 7, 073509-1--073509-12  (2005).
%
\bibitem{IMQ01}Kawamura, K.: 
Indefinite-metric quantum field theory and operator algebra,
math.OA/0608076 (2006).
%
\bibitem{RBS01}Kawamura, K.:
Recursive boson system in the Cuntz algebra ${\cal O}_{\infty}$,
J. Math. Phys.  {\bf 48}, No. 9,  093510-1--093510-16 (2007).
%
\bibitem{TS01}Kawamura, K.:
A tensor product of representations of Cuntz algebras,
Lett. Math. Phys., to appear.
%
\bibitem{KelVau}Kelley, J. L.  and Vaught, R. L.: 
The positive cone in Banach algebras,
 Trans. Amer. Math. Soc. {\bf 74} (1953), 44--55.
%
\bibitem{MMSV}Mnatskanova, M., Morchio,  G., Strocchi,  F.  and Vernov, Yu.: 
Representations of CCR algebras in Krein spaces of entire functions,
Lett. Math. Phys. {\bf 65} (2003) 159--172.
%
%
\bibitem{MS1980}
Morchio, G. and Strocchi, F.: 
Infrared singularities, vacuum structure and pure phases 
in local quantum field theory,
Ann. Inst. H. Poincare Sect. A (N.S.) {\bf 33} (1980), no. 3, 251--282.
%
\bibitem{MPS1990}Morchio, G., Pierotti, D. and Strocchi, F.: 
Infrared and vacuum structure in two-dimensional models of
local quantum field theory. The massless scalar field,
J. Math. Phys. {\bf 31} (1990) 1467--1477.
%
\bibitem{MPS1992}Morchio, G., Pierotti, D. and Strocchi, F.: 
Infrared and vacuum structure in two-dimensional models of
local quantum field theory. II. Fermion bosonization,
J. Math. Phys. {\bf 33} (1992) 777--790.
%
\bibitem{Nagy}Nagy,  K. L.: 
{\it State vector spaces with indefinite metric in quantum field theory}, 
P. Noordhoff Ltd., Groningen; Akademiai Kiado, Budapest 1966.
%
\bibitem{Naimark}Na\u\i mark, M. A.: 
Rings with involutions, 
Uspehi Matem. Nauk (N.S.) 3, (1948). no. 5(27), 52--145.
%
\bibitem{NakagamiTomita}Nakagami, Y.  and Tomita, M.: 
Triangular matrix representation for self-adjoint operators in Krein spaces,
Japan. J. Math. {\bf 14}, No 1 (1988) 165--202.
%
\bibitem{Nakanishi}Nakanishi, N.: 
Indefinite-metric quantum field theory,
Suppl. Prog. Theor. Phys. {\bf 51} (1972) 1--95.
%
\bibitem{Pauli}Pauli, W.: 
On Dirac's new method of field quantization,
Rev. Modern. Phys. {\bf 15} (1943) 175--207.
%
\bibitem{Ped}Pedersen,  G. K.: 
{\it $C^*$-algebras and their automorphism groups},
Academic Press Incorporated (1979).
%
\bibitem{Pontrjagin}Pontrjagin, L. S.: 
Hermitian operators in spaces with indefinite metric,
Izv. Akad.  Nauk SSSR Ser. Mat. {\bf 8} (1944) 243--280.
%
\bibitem{Sobolev}Sobolev, S. L.: 
The motion of a symmetric top containing a cavity filled 
with a liquid,
$\check{Z}$. Prikl. Meh. i Tehn. Fiz. (1960), no 3, 20--55.
%
\bibitem{Solombrino}Solombrino, L.: 
Weak pseudo-Hermiticity and antilinear commutant,
J. Math. Phys. Vol {\bf 43},  11 (2002) 5439--5445.
%
\bibitem{Strocchi}Strocchi, F.: 
{\it Elements of quantum mechanics of infinite systems,} 
International School for Advanced Studies Lecture Series, 3. 
World Scientific Publishing Co., Singapore, 1985.  
%
\bibitem{Thaller}Thaller, B.: 
{\it The Dirac equation,}
Springer-Verlag, 1992.
%
\bibitem{Tomita}Tomita, M.: 
Operators and operator algebras in Krein spaces, I--Spectral 
analysis in Pontrjagin space--,
RIMS-k\=oky\=uroku, {\bf 398} (1980), 131--158.
%
\bibitem{Weyl}Weyl,  H.: 
Theorie der Darstellung kontinuierlicher halb-einfacher 
Gruppen durch lineare Transformationen. I, 
Math. Z. {\bf 23} (1925), no. 1, 271--309.
%
\bibitem{Weyl02}Weyl, H.: 
Theorie der Darstellung kontinuierlicher halb-einfacher Gruppen 
durch lineare Transformationen. II,
 Math. Z. {\bf 24} (1926), no. 1, 328--376.
%
\bibitem{Weyl03}Weyl, H.: 
Theorie der Darstellung kontinuierlicher halb-einfacher 
Gruppen durch lineare Transformationen. III, 
Math. Z. {\bf 24} (1926), no. 1, 377--395. 
\end{thebibliography}
\end{document}